\documentclass[conference]{IEEEtran}
\IEEEoverridecommandlockouts
\usepackage{cite}
\usepackage{amsmath,amssymb,amsfonts}
\usepackage{algorithmic}
\usepackage{graphicx}
\usepackage{textcomp}
\usepackage{xcolor}
\def\BibTeX{{\rm B\kern-.05em{\sc i\kern-.025em b}\kern-.08em
    T\kern-.1667em\lower.7ex\hbox{E}\kern-.125emX}}

\begin{document}

\title{BESS Participation Planning for Provision of Grid Services in Energy and Regulation Markets}

\author{\IEEEauthorblockN{Zeenat Hameed, Chresten Træholt}
\IEEEauthorblockA{\textit{Technical University of Denmark, Department of Wind and Energy Systems}\\
}}

\maketitle

\begin{abstract}
This paper presents a stochastic optimization model for planning the participation of battery energy storage systems (BESSs) in energy and regulation markets. The proposed model quantifies and compares the business value of single and multi-market BESS services by accounting for their bid submission and acceptance procedures, pricing mechanisms and revenue streams, and penalty payments and degradation costs. By modelling the price uncertainties using market price scenarios (MPSs) and considering representative days for different price seasonalities, the model outputs the percentage of hours each market service can be targeted to maximize profits. This helps in quantifying the impact of operational and market requirements of different services on choice of markets for BESSs. It also helps in determining the approximate costs and revenues that may be accrued by the BESS owners by choosing combinations of available services under different price conditions. The model thus overcomes the key limitations of previous studies that were mainly conducted from a controller design viewpoint and were thus more focused on the operational control of BESSs. The proposed model is generalizable and extendable to BESS service provision in multiple markets of different regions.
\end{abstract}


\section{Introduction}

Environmental concerns and climate‐change awareness have accelerated integration of renewable energy sources (RES) into power systems \cite{b1}. However, RES intermittency and low synchronous inertia undermine grid stability and reliability. Battery energy storage systems (BESSs) alleviate these issues—especially in regulation markets—by providing fast response and flexible control for frequency support and ancillary services \cite{b2,b3,b4,b5,b6}. Traditionally, system operators sourced regulation from fossil plants at administered rates \cite{b7}, but BESSs cannot hold capacity indefinitely due to limited state‐of‐charge \cite{b8}. Although capital costs are declining, BESS payback periods remain long, making optimal service scheduling essential for profitability \cite{b9}.

Most existing work optimizes BESS operation from a controller‐design perspective, treating BESSs as standalone multi‐service units. For example, \cite{b10} proposes a control strategy for simultaneous energy arbitrage and ancillary services; \cite{b11} minimizes a commercial consumer’s bill by combining peak shaving with frequency regulation; \cite{b12} integrates day‐ahead dispatch with frequency and voltage regulation; and \cite{b13} replaces droop control with output‐regulation theory for voltage/frequency support. Although these approaches improve operational algorithms and correct forecast errors, they generally omit market‐specific parameters—bid submission and acceptance rules, pricing mechanisms, penalties for non‐delivery, and diverse revenue streams—which are critical to assessing BESS value from a business perspective \cite{b14,b15}.

Furthermore, most literature focuses on U.S.\ markets with locational marginal pricing and pay‐for‐performance schemes. European markets, by contrast, employ zonal pricing and pay for availability rather than performance \cite{b16,b17, b18, b19}. To address this gap, we develop a stochastic optimization model that incorporates bid procedures, pricing rules, penalty structures, and degradation costs to quantify the business value of BESS participation across multiple electricity markets.

\section{Nordic Market Mechanisms}

The Nordic market consists of a day‐ahead spot market, an intra‐day market (both administered by Nord Pool Spot), a regulation market, and a balancing market (operated by national TSOs: Svenska Kraftnät, Statnett, Fingrid, Energinet). Denmark uses zonal pricing, in line with other European regions \cite{b1,b2,b3}.

\paragraph{Spot and Intra‐day}  
Participants submit hourly bids 24 h ahead; intra‐day adjustments are allowed until one hour before delivery. Both use a pay‐as‐clear mechanism: the marginal clearing price applies per MWh.

\paragraph{Regulation}  
Ancillary services—FCR-N, FCR-D, aFRR, mFRR, FFR—are procured via voluntary bids 1–2 days ahead. TSOs activate bids based on system needs:
\begin{itemize}
  \item \textbf{FCR-N:} Symmetrical response when \(f\le49.9\)\,Hz or \(f\ge50.1\)\,Hz.
  \item \textbf{FCR-D:} Separate up/down bids for \(f\le49.5\)\,Hz or \(f\ge50.5\)\,Hz.
\end{itemize}
BESSs target FCR-N and FCR-D due to fast ramping. Bids specify MW capacity and are paid pay‐as‐bid for availability; activated energy settles at balancing prices (aFRR/mFRR). Failure to deliver incurs penalties equal to lost availability plus imbalance charges \cite{b4,b5}.

\begin{figure}[ht]
    \centering
    % Placeholder for the actual binary decision diagram
    \includegraphics[width=0.8\linewidth]{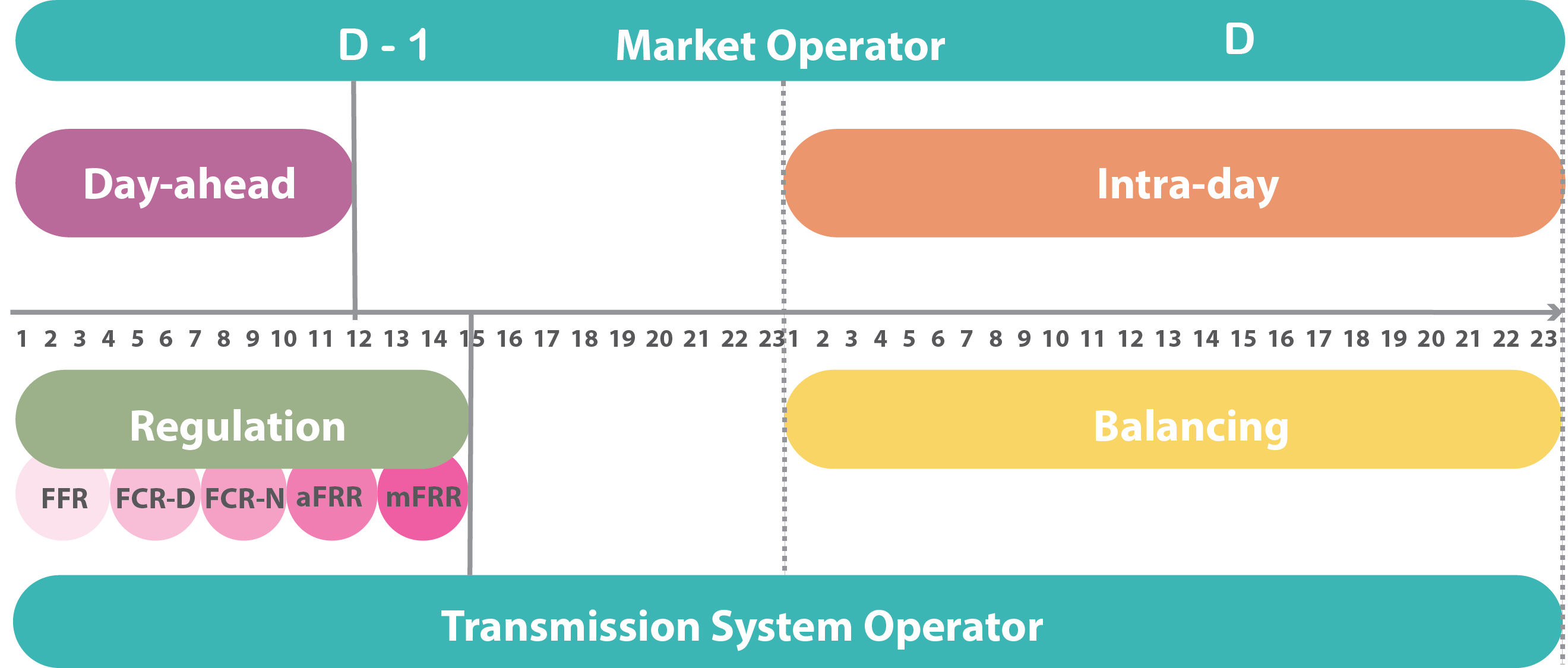}
    \caption{Nordic markets\label{fig:markets_2}}
\end{figure}
\paragraph{Balancing}  
Deviations from spot/intra‐day schedules are settled via single‐ or double‐price mechanisms. If a regulation bid is accepted but not fulfilled, imbalance prices apply, potentially exceeding initial revenues \cite{b6}.

\paragraph{Payment Structures}  
\begin{itemize}
  \item \textbf{Spot (S-DCH/S-CH):} Pay-as-clear per MWh—discharge (S-DCH) earns the spot price; charge (S-CH) pays it.
  \item \textbf{FCR-N/FCR-D:} Pay-as-bid availability per MW; if activated, net discharge settles at up-regulation price, net charge at down-regulation price. Bids are in one-hour blocks; penalties cover lost availability and imbalance charges.
\end{itemize}

This summary informs our model’s bid logic, revenue, and penalty formulations.

\section{Method}

\subsection{Droop Control}
In the Nordic system, the day‐ahead spot market handles energy trades for the next 24 h, while real‐time balance relies on regulation and balancing markets. TSOs activate FCR‐N and FCR‐D bids (submitted 1–2 days ahead) for frequency support; balancing penalizes deviations from declared schedules. Although frequency varies millisecond‐wise, we sample per minute. BESS participation in FCR‐N, FCR‐D, and day‐ahead spot is modeled via per‐minute energy requirements \(E_{stm}^{\mathrm{dch}}\) (upward) and \(E_{stm}^{\mathrm{ch}}\) (downward) for each scenario \(s\in S\), period \(t\in T\), and market \(m\in M\). These are derived from droop control points \((P_{m1},P_{\max}),\,(P_{m2},0),\,(P_{m3},0),\,(P_{m4},-P_{\max})\):
\[
(P_{m1}, P_{\max}),\; (P_{m2}, 0),\; (P_{m3}, 0),\; (P_{m4}, -P_{\max}).
\]
For \(m=N\) (FCR‐N): \(P_{m1}=49.9\), \(P_{m2}=P_{m3}=50.0\), \(P_{m4}=50.1\). For \(m=D\) (FCR‐D): \(P_{m1}=49.5\), \(P_{m2}=49.9\), \(P_{m3}=50.1\), \(P_{m4}=50.5\). Applying linear interpolation,
\begin{equation}
E_{stm}^{\mathrm{dch}} \;=\; \frac{1}{60}\,
\begin{cases}
0, & f_{st}^m < P_{m1},\\
P_{\max} - \dfrac{P_{\max}}{P_{m2}-P_{m1}}\,(f_{st}^m - P_{m1}), & P_{m1}\le f_{st}^m\le P_{m2},\\
0, & \text{otherwise},
\end{cases}
\end{equation}
\begin{equation}
E_{stm}^{\mathrm{ch}} \;=\; \frac{1}{60}\,
\begin{cases}
0, & f_{st}^m < P_{m3},\\
\dfrac{P_{\max}}{P_{m4}-P_{m3}}\,\bigl(f_{st}^m - P_{m3}\bigr), & P_{m3}\le f_{st}^m\le P_{m4},\\
0, & f_{st}^m > P_{m4}.
\end{cases}
\end{equation}
Since only FCR‐D upward is considered, \(E_{stm}^{\mathrm{ch}}=0\) for \(m=D\).

For spot markets, with fixed bid energy \(E_{\text{spot}}\):
\begin{align}
&\text{S‐DCH (sell):}\quad E_{stm}^{\mathrm{dch}} = \frac{E_{\text{spot}}}{60}, \quad E_{stm}^{\mathrm{ch}} = 0,\\
&\text{S‐CH (buy):}\quad E_{stm}^{\mathrm{ch}} = \frac{E_{\text{spot}}}{60}, \quad E_{stm}^{\mathrm{dch}} = 0.
\end{align}

These \(E_{stm}^{\mathrm{dch}}\) and \(E_{stm}^{\mathrm{ch}}\) values feed into the optimizer.
\begin{figure}[ht]
    \centering
    % Placeholder for the actual binary decision diagram
    \includegraphics[width=0.8\linewidth]{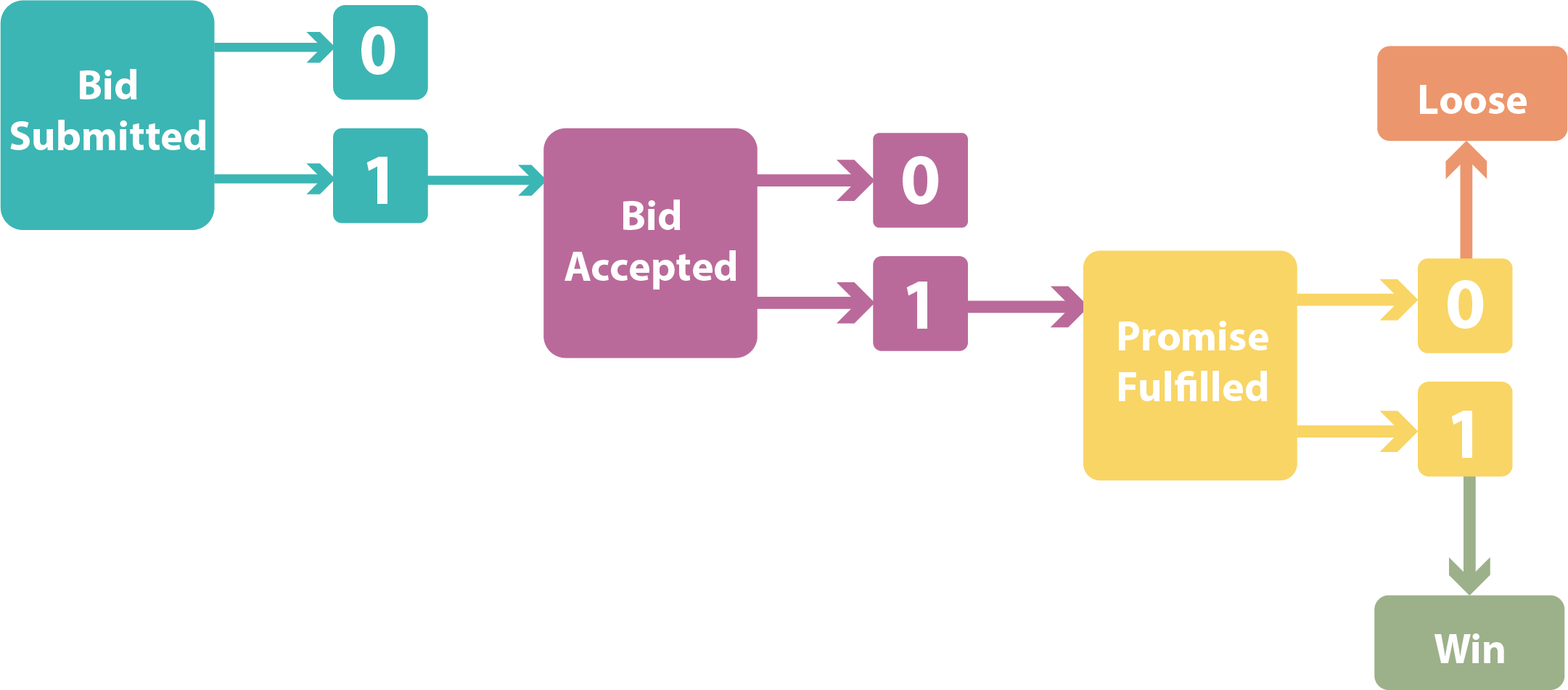}
    \caption{Binary decisions\label{fig:binary}}
\end{figure}
\subsection{Binary Decisions}

At each hour \(h\in H\), for each market \(m\in M\) and scenario \(s\in S\), the optimizer makes three binary decisions (Fig.~\ref{fig:binary}):

1. \(\,x_{hm}^{\mathrm{bid}}\in\{0,1\}\): whether to submit a bid. This depends on the BESS state‐of‐charge \(z_{st}^{\mathrm{SOC}}\) and market requirements \(E_{stm}^{\mathrm{dch}}\), \(E_{stm}^{\mathrm{ch}}\).

2. \(\,x_{shm}^{\mathrm{bid,acc}}\in\{0,1\}\): whether the bid is accepted. If \(x_{hm}^{\mathrm{bid}}=1\), acceptance requires the bid price \(x_{h}^{\mathrm{price}}\) to meet the scenario‐specific threshold \(P_{shm}^{\max}\).

3. \(\,w_{sh}^{\mathrm{ok}}\in\{0,1\}\): whether the BESS fulfills its energy promise. If \(x_{shm}^{\mathrm{bid,acc}}=1\), this depends on matching per‐minute dispatch \(z_{stm}^{\mathrm{dch}}\), \(z_{stm}^{\mathrm{ch}}\) to the market’s \(E_{stm}^{\mathrm{dch}}\), \(E_{stm}^{\mathrm{ch}}\); failure incurs a penalty.

\subsection{Regulation Market Payments}

Since availability payments of FCR-N and FCR-D markets ($m \in M_{\text{freq}}$) operate on pay-as-bid mechanism and are applied per MW, $x_{h}^{\text{price}}$ is multiplied by the bid power ($P_{\max}$). To simplify, we assume the BESS always bids with $P_{\max} = 0.9$~MW. Moreover, they also depend on binary decision variable $x_{shm}^{\text{bid,acc}}$ as discussed in the previous section. Therefore, the payment for scenario $s$, hour $h$ is given by:
\begin{equation}
    P_{\max} \, x_{h}^{\text{price}} \sum_{m \in M_{\text{freq}}} x_{shm}^{\text{bid,acc}}.
\end{equation}

Since this expression is not linear, we use the Big-$M$ method to linearize it. We introduce $w_{sh}^{\text{avail}} \ge 0$ and let:
\begin{align}
    w_{sh}^{\text{avail}} &= P_{\max} \, x_{h}^{\text{price}} \sum_{m \in M_{\text{freq}}} x_{shm}^{\text{bid,acc}}, \\
    w_{sh}^{\text{avail}} &\le P_{\max} \, x_{h}^{\text{price}}, \\
    w_{sh}^{\text{avail}} &\le M \sum_{m \in M_{\text{freq}}} x_{shm}^{\text{bid,acc}}, \\
    w_{sh}^{\text{avail}} &\ge P_{\max} \, x_{h}^{\text{price}} - M\left(1 - \sum_{m \in M_{\text{freq}}} x_{shm}^{\text{bid,acc}}\right).
\end{align}

\subsection{Spot Market Payments}

Since spot market transactions of S-DCH and S-CH operate on pay-as-clear mechanism and are applied per MWh, $P_{shm}^{\max}$ is multiplied by the bid energy ($E_{\text{spot}}$). To simplify, we assume the BESS always bids with $E_{\text{spot}} = 0.4$~MWh. Moreover, they also depend on binary decision variable $x_{shm}^{\text{bid,acc}}$. Therefore, the payment for scenario $s$, hour $h$ in S-DCH is given by:
\begin{equation}
    E_{\text{spot}} \, P_{shm}^{\max} \, x_{sh,\text{S-DCH}}^{\text{bid,acc}}.
\end{equation}
To linearize, we introduce $w_{sh}^{\text{spot,dch}} \ge 0$ and let:
\begin{align}
    w_{sh}^{\text{spot,dch}} &= E_{\text{spot}} \, P_{shm}^{\max} \, x_{sh,\text{S-DCH}}^{\text{bid,acc}}, \\
    w_{sh}^{\text{spot,dch}} &\le E_{\text{spot}} \, P_{shm}^{\max}, \\
    w_{sh}^{\text{spot,dch}} &\le M \, x_{sh,\text{S-DCH}}^{\text{bid,acc}}, \\
    w_{sh}^{\text{spot,dch}} &\ge E_{\text{spot}} \, P_{shm}^{\max} - M\left(1 - x_{sh,\text{S-DCH}}^{\text{bid,acc}}\right).
\end{align}

On the other hand, the amount we have to pay for scenario $s$, hour $h$ in S-CH is given by:
\begin{equation}
    E_{\text{spot}} \, P_{shm}^{\max} \, x_{sh,\text{S-CH}}^{\text{bid,acc}}.
\end{equation}
It is linearized by introducing $w_{sh}^{\text{spot,ch}} \ge 0$ and letting:
\begin{align}
    w_{sh}^{\text{spot,ch}} &= E_{\text{spot}} \, P_{shm}^{\max} \, x_{sh,\text{S-CH}}^{\text{bid,acc}}, \\
    w_{sh}^{\text{spot,ch}} &\le E_{\text{spot}} \, P_{shm}^{\max}, \\
    w_{sh}^{\text{spot,ch}} &\le M \, x_{sh,\text{S-CH}}^{\text{bid,acc}}, \\
    w_{sh}^{\text{spot,ch}} &\ge E_{\text{spot}} \, P_{shm}^{\max} - M\left(1 - x_{sh,\text{S-CH}}^{\text{bid,acc}}\right).
\end{align}

\subsection{Penalty Payments}

When \(x_{shm}^{\mathrm{bid,acc}}=1\), two cases arise depending on \(w_{sh}^{\mathrm{ok}}\in\{0,1\}\):

- If \(w_{sh}^{\mathrm{ok}}=1\), the BESS earns:
  \[
    w_{sh}^{\mathrm{avail}} \;+\; w_{sh}^{\mathrm{spot,dch}} \;-\; w_{sh}^{\mathrm{spot,ch}}.
  \]
- If \(w_{sh}^{\mathrm{ok}}=0\), it incurs penalties. For \(m\in M_{\mathrm{freq}}\), penalty is
  \[
    w_{sh}^{\mathrm{avail}},
  \]
  and for spot markets:
  \[
    \text{S‐DCH: }C_{sh}^{\mathrm{up}}\,E_{\mathrm{spot}}, 
    \qquad
    \text{S‐CH: }C_{sh}^{\mathrm{down}}\,E_{\mathrm{spot}}.
  \]
  (We apply full \(E_{\mathrm{spot}}\) for simplicity.)

Hence the net revenue in scenario \(s\), hour \(h\) is
\[
  \begin{cases}
    w_{sh}^{\mathrm{avail}} + w_{sh}^{\mathrm{spot,dch}} - w_{sh}^{\mathrm{spot,ch}}, 
      & w_{sh}^{\mathrm{ok}} = 1, \\[0.5ex]
    -\,\bigl[w_{sh}^{\mathrm{avail}} + C_{sh}^{\mathrm{up}}\,E_{\mathrm{spot}}
      + C_{sh}^{\mathrm{down}}\,E_{\mathrm{spot}}\bigr], 
      & w_{sh}^{\mathrm{ok}} = 0.
  \end{cases}
\]

To linearize products with \(w_{sh}^{\mathrm{ok}}\), define auxiliary variables:

\paragraph{Availability “win” (when \(m\in M_{\mathrm{freq}}\))}  
\[
  w_{sh}^{\mathrm{won,avail}} = w_{sh}^{\mathrm{avail}} \,w_{sh}^{\mathrm{ok}},
\]
subject to
\begin{align*}
  w_{sh}^{\mathrm{won,avail}} &\le w_{sh}^{\mathrm{avail}},\\
  w_{sh}^{\mathrm{won,avail}} &\le M\,w_{sh}^{\mathrm{ok}},\\
  w_{sh}^{\mathrm{won,avail}} &\ge w_{sh}^{\mathrm{avail}} - M\,(1 - w_{sh}^{\mathrm{ok}}).
\end{align*}

\paragraph{Availability “loss” (when \(m\in M_{\mathrm{freq}}\))}  
\[
  w_{sh}^{\mathrm{lost,avail}} = w_{sh}^{\mathrm{avail}}\,(1 - w_{sh}^{\mathrm{ok}}),
\]
subject to
\begin{align*}
  w_{sh}^{\mathrm{lost,avail}} &\le w_{sh}^{\mathrm{avail}},\\
  w_{sh}^{\mathrm{lost,avail}} &\le M\,(1 - w_{sh}^{\mathrm{ok}}),\\
  w_{sh}^{\mathrm{lost,avail}} &\ge w_{sh}^{\mathrm{avail}} - M\,w_{sh}^{\mathrm{ok}}.
\end{align*}

\paragraph{Spot “win” for discharge (when \(m\in M_{\mathrm{S-DCH}}\))}  
\[
  w_{sh}^{\mathrm{won,spot,dch}} = w_{sh}^{\mathrm{spot,dch}}\,w_{sh}^{\mathrm{ok}},
\]
subject to
\begin{align*}
  w_{sh}^{\mathrm{won,spot,dch}} &\le w_{sh}^{\mathrm{spot,dch}},\\
  w_{sh}^{\mathrm{won,spot,dch}} &\le M\,w_{sh}^{\mathrm{ok}},\\
  w_{sh}^{\mathrm{won,spot,dch}} &\ge w_{sh}^{\mathrm{spot,dch}} - M\,(1 - w_{sh}^{\mathrm{ok}}).
\end{align*}

\paragraph{Spot “loss” for discharge (when \(m\in M_{\mathrm{S-DCH}}\))}  
\[
  w_{sh}^{\mathrm{lost,spot,dch}} = C_{sh}^{\mathrm{up}}\,E_{\mathrm{spot}}\,(1 - w_{sh}^{\mathrm{ok}}),
\]
subject to
\begin{align*}
  w_{sh}^{\mathrm{lost,spot,dch}} &\le C_{sh}^{\mathrm{up}}\,E_{\mathrm{spot}},\\
  w_{sh}^{\mathrm{lost,spot,dch}} &\le M\,(1 - w_{sh}^{\mathrm{ok}}),\\
  w_{sh}^{\mathrm{lost,spot,dch}} &\ge C_{sh}^{\mathrm{up}}\,E_{\mathrm{spot}} - M\,w_{sh}^{\mathrm{ok}}.
\end{align*}

\paragraph{Spot “loss” for charge (when \(m\in M_{\mathrm{S-CH}}\))}  
\[
  w_{sh}^{\mathrm{lost,spot,ch}} = C_{sh}^{\mathrm{down}}\,E_{\mathrm{spot}}\,(1 - w_{sh}^{\mathrm{ok}}),
\]
subject to
\begin{align*}
  w_{sh}^{\mathrm{lost,spot,ch}} &\le C_{sh}^{\mathrm{down}}\,E_{\mathrm{spot}},\\
  w_{sh}^{\mathrm{lost,spot,ch}} &\le M\,(1 - w_{sh}^{\mathrm{ok}}),\\
  w_{sh}^{\mathrm{lost,spot,ch}} &\ge C_{sh}^{\mathrm{down}}\,E_{\mathrm{spot}} - M\,w_{sh}^{\mathrm{ok}}.
\end{align*}

\subsection{Objective Function and Constraints}

The objective function maximizes the total expected revenue comprising of earnings and penalties of each market, and the degradation costs. It is given as:
\begin{equation}
\begin{aligned}
\max\;\sum_{s \in S} p_s \Biggl[ 
    &\sum_{h \in H} \bigl(w_{sh}^{\mathrm{won,avail}} 
        +\,w_{sh}^{\mathrm{won,s\text{-}dch}}\bigr) \\[0.5ex]
    &\quad - \sum_{h \in H} \bigl(w_{sh}^{\mathrm{lost}} 
        +\,w_{sh}^{\mathrm{lost,s\text{-}dch}} 
        +\,w_{sh}^{\mathrm{lost,s\text{-}ch}}\bigr) \\[0.5ex]
    &\quad + \sum_{h \in H} w_{sh}^{\mathrm{energy}} \\[0.5ex]
    &\quad - \sum_{t \in T} \sum_{m \in M} 
            C_{shtm}^{\mathrm{degradation}}\;x_{shm}^{\mathrm{bid,acc}}
\Biggr]
\end{aligned}
\end{equation}

\noindent
\resizebox{\columnwidth}{!}{%
  \begin{minipage}{\columnwidth}
    \scriptsize
    % (Removed \setlength{\mathindent}{0pt}, since it’s undefined here.)
    \allowdisplaybreaks
    \begin{alignat*}{2}
      &\sum_{m \in M} x_{hm}^{\mathrm{bid}}
        \;\le\; 1,
        && \quad \forall\,h \in H, \\[0.5ex]
      &x_{shm}^{\mathrm{bid,acc}}
        \;\le\; x_{hm}^{\mathrm{bid}},
        && \quad \forall\,s \in S,\,h \in H,\,m \in M, \\[0.5ex]
      &\sum_{m \in M}\bigl(P_{shm}^{\max}\,x_{hm}^{\mathrm{bid}}\bigr)
        \;-\; x_{h}^{\mathrm{price}}
        \;\le\; M\,\sum_{m \in M} x_{shm}^{\mathrm{bid,acc}},
        && \quad \forall\,s \in S,\,h \in H, \\[0.5ex]
      &x_{h}^{\mathrm{price}} \;-\; P_{shm}^{\max}
        \;\le\; M\,\bigl(1 - x_{shm}^{\mathrm{bid,acc}}\bigr),
        && \quad \forall\,s \in S,\,h \in H,\,m \in M, \\[0.5ex]
      &x_{h}^{\mathrm{price}}
        \;\le\; \mathrm{Bid}_{\max},
        && \quad \forall\,h \in H, \\[0.5ex]
      &x_{h}^{\mathrm{price}}
        \;\ge\; \mathrm{Bid}_{\min}\,\sum_{m \in M} x_{hm}^{\mathrm{bid}},
        && \quad \forall\,h \in H, \\[1ex]
      &z_{stm}^{\mathrm{dch}}
        \;=\; 0,
        && \quad \forall\,s \in S,\,t \in T,\,m \in M:\;P_{m2}<f_{st}^{m}, \\[0.5ex]
      &z_{stm}^{\mathrm{ch}}
        \;=\; 0,
        && \quad \forall\,s \in S,\,t \in T,\,m \in M:\;P_{m2}<f_{st}^{m}, \\[0.5ex]
      &z_{stm}^{\pm}
        \;=\; 0,
        && \quad \forall\,s \in S,\,t \in T,\,m \in M:\;P_{m2}<f_{st}^{m}, \\[1ex]
      &z_{stm}^{\mathrm{dch}} + s_{stm}^{\mathrm{dch}}
        \;=\; E_{stm}^{\mathrm{dch}}\,x_{shm}^{\mathrm{bid,acc}},
        && \quad \forall\,s \in S,\,t \in T,\,m \in M, \\[0.5ex]
      &z_{stm}^{\mathrm{ch}} + s_{stm}^{\mathrm{ch}}
        \;=\; E_{stm}^{\mathrm{ch}}\,x_{shm}^{\mathrm{bid,acc}},
        && \quad \forall\,s \in S,\,t \in T,\,m \in M, \\[0.5ex]
      &z_{stm}^{\pm}
        \;=\; z_{stm}^{\mathrm{dch}} \;-\; z_{stm}^{\mathrm{ch}},
        && \quad \forall\,s \in S,\,t \in T,\,m \in M, \\[1ex]
      &z_{s0}^{\mathrm{SOC}}
        \;=\; M_{0},
        && \quad \forall\,s \in S, \\[0.5ex]
      &E_{\min} \;\le\; z_{st}^{\mathrm{SOC}} 
          \;-\; \sum_{m \in M} z_{stm}^{\pm}
        \;\le\; E_{\max},
        && \quad \forall\,s \in S,\,t \in T, \\[0.5ex]
      &z_{st}^{\mathrm{SOC}}
        \;=\; z_{s,t-1}^{\mathrm{SOC}} 
              \;-\; \sum_{m \in M} z_{s,t-1,m}^{\pm},
        && \quad \forall\,s \in S,\,t \in T, \\[1ex]
      &\sum_{m \in M} \sum_{t : h_{t} = h} \bigl(s_{stm}^{\mathrm{dch}} + s_{stm}^{\mathrm{ch}}\bigr)
        \;\le\; M\,\bigl(1 - w_{sh}^{\mathrm{ok}}\bigr),
        && \quad \forall\,s \in S,\,h \in H, \\[1ex]
      &x_{hm}^{\mathrm{bid}} 
        \;\in\; \{0,1\}, 
        && \quad \forall\,h \in H,\,m \in M, \\[0.5ex]
      &x_{h}^{\mathrm{price}} 
        \;\ge\; 0, 
        && \quad \forall\,h \in H, \\[0.5ex]
      &x_{shm}^{\mathrm{bid,acc}} 
        \;\in\; \{0,1\}, 
        && \quad \forall\,s \in S,\,h \in H,\,m \in M, \\[0.5ex]
      &z_{st}^{\mathrm{SOC}} 
        \;\ge\; 0, 
        && \quad \forall\,s \in S,\,t \in T, \\[0.5ex]
      &z_{stm}^{\mathrm{dch}} 
        \;\ge\; 0, 
        && \quad \forall\,s \in S,\,t \in T,\,m \in M, \\[0.5ex]
      &z_{stm}^{\mathrm{ch}} 
        \;\ge\; 0, 
        && \quad \forall\,s \in S,\,t \in T,\,m \in M, \\[0.5ex]
      &z_{stm}^{\pm} 
        \;\in\; \mathbb{R}, 
        && \quad \forall\,s \in S,\,t \in T,\,m \in M, \\[0.5ex]
      &s_{stm}^{\mathrm{dch}} \;\ge\; 0,\; s_{stm}^{\mathrm{ch}} \;\ge\; 0
        &,\quad &\forall\,s \in S,\,t \in T,\,m \in M, \\[1ex]
      &w_{sh}^{\mathrm{avail}} \;\ge\; 0,\; w_{sh}^{\mathrm{spot,dch}} \;\ge\; 0,\; w_{sh}^{\mathrm{spot,ch}} \;\ge\; 0
        &,\quad &\forall\,s \in S,\,h \in H, \\[0.5ex]
      &w_{sh}^{\mathrm{ok}}
        \;\in\; \{0,1\}, 
        && \quad \forall\,s \in S,\,h \in H, \\[0.5ex]
      &w_{sh}^{\mathrm{won,avail}} \;\ge\; 0,\; w_{sh}^{\mathrm{lost,avail}} \;\ge\; 0
        &,\quad &\forall\,s \in S,\,h \in H, \\[0.5ex]
      &w_{sh}^{\mathrm{won,spot,dch}} \;\ge\; 0,\; w_{sh}^{\mathrm{lost,spot,dch}} \;\ge\; 0,\; w_{sh}^{\mathrm{lost,spot,ch}} \;\ge\; 0
        &,\quad &\forall\,s \in S,\,h \in H, \\[0.5ex]
      &w_{sh}^{\mathrm{energy}}
        \;\ge\; 0, 
        && \quad \forall\,s \in S,\,h \in H, \\[0.5ex]
      &s_{st}^{\mathrm{SOC\text{-}target},+},\, s_{st}^{\mathrm{SOC\text{-}target},-} \;\ge\; 0
        &,\quad &\forall\,s \in S,\,t \in \{\,h_{N} : h \in \{1,\dots,23\}\}. 
    \end{alignat*}
  \end{minipage}%
}% end \resizebox

\section{Results and Discussion}

To assess BESS performance, historical data of one representative week of 2024 ($W_{\text{rep}}^{2024}$) is considered. Each day of $W_{\text{rep}}^{2024}$ is considered a scenario that consists of per-hour market prices, and per-minute grid-frequencies. The optimizer thus approximates optimal daily BESS performance in $W_{\text{rep}}^{2024}$.

\subsection{BESS Services in Single Markets}

First, single-market BESS participation is considered. Applying DCS from (1)–(7) to $f_{st}$ of $W_{\text{rep}}^{2024}$ gives $E_{stm}^{dch}$. The summary of the observed parameters is shown in Fig.~\ref{fig:single_bid}. Both cases when degradation costs are not considered (no-DC) shown in aqua, and when degradation costs are considered (DC) shown in orange are represented.

When $m = N$, for no-DC case, bids are submitted in 100\% of the hours and accepted when $x_{h}^{\text{price}} \le P_{s3hN}^{\max}$. This results in the BESS being in FCR-N market 92\% of the time, which is also illustrated in Fig.~\ref{fig:participation_single}. 8\% hours are idle. Contrarily, for DC case, both bid-submission and acceptance hours are reduced to 79\%, and 38\% to avoid charge/discharge cycles that incur higher degradation costs. However, relative to no-DC case, even though BESS participation in FCR-N market is reduced by 34\%, its profits are only reduced by 9.6\% as evident in Fig.~\ref{fig:participation_single}. One reason for this is the avoidance of penalty costs in the DC case that were mainly incurred during high degradation cost hours.

When $m = D$, both in DC and no-DC case, bids are submitted and accepted 100\% of the time. The reason is the low $E_{stD}^{dch}$ requirements of FCR-D market, which translate to significantly low charge/discharge cycles thereby resulting in negligible degradation costs. Even with higher degradation costs of the latter, its overall profits still remain 32\% higher than the former. When $m = \text{S-DCH} + \text{S-CH}$, for no-DC case, bids are respectively submitted and accepted 54\%, and 12\% of the time, which drops even further to 38\% and 4\% in the DC case. Moreover, in the DC case, the BESS only participates in S-DCH and avoids S-CH altogether.

\begin{figure}[ht]
    \centering
    % Placeholder for bid/submission figure
    \includegraphics[width=\linewidth]{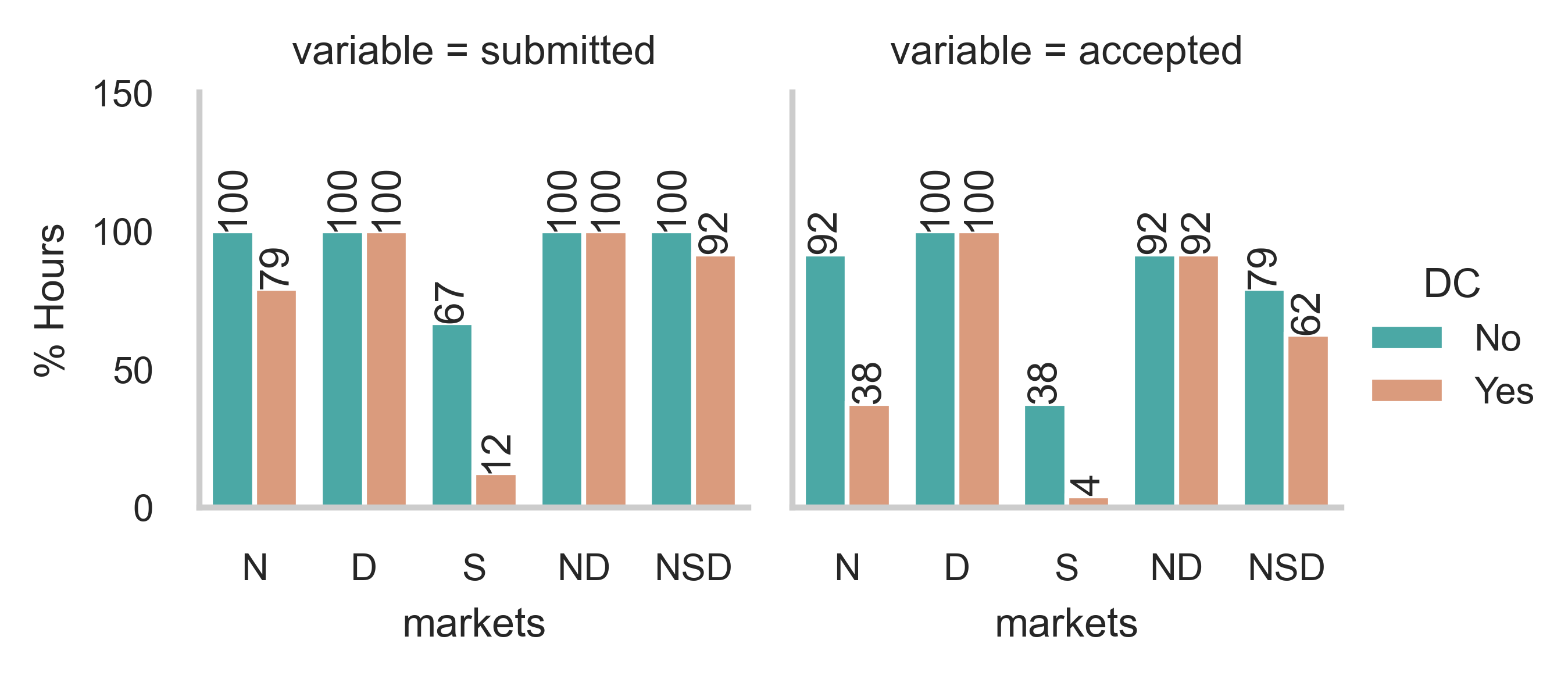}
    \caption{Bid submitted and accepted hours for single-market services (scenario $s_3$)\label{fig:opt_r11}}
\end{figure}

\begin{figure}[ht]
    \centering
    % Placeholder for costs and earnings figure
    \includegraphics[width=\linewidth]{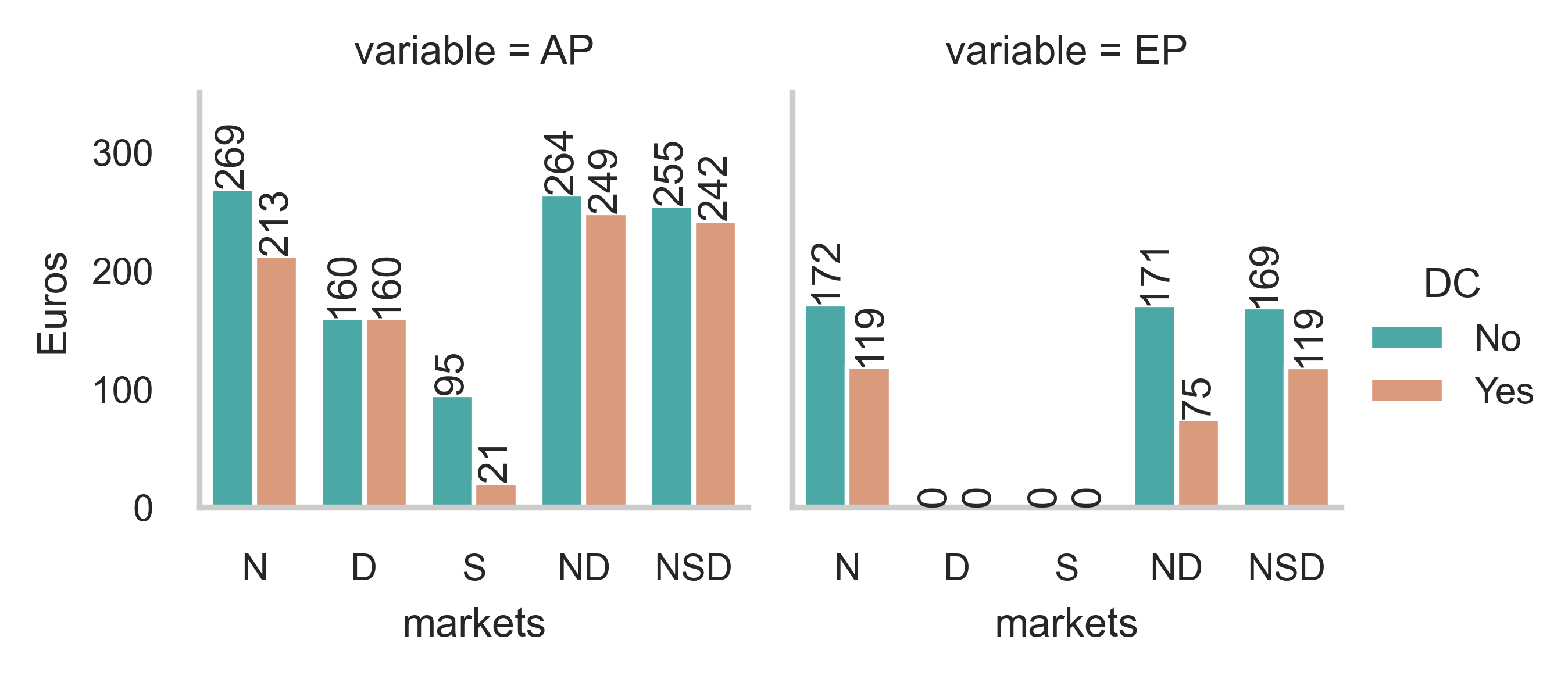}
    \caption{Costs and earnings for single-market participation (scenario $s_3$)\label{fig:opt_r12}}
\end{figure}
\begin{figure}[ht]
    \centering
    % Placeholder for costs and earnings figure
    \includegraphics[width=\linewidth]{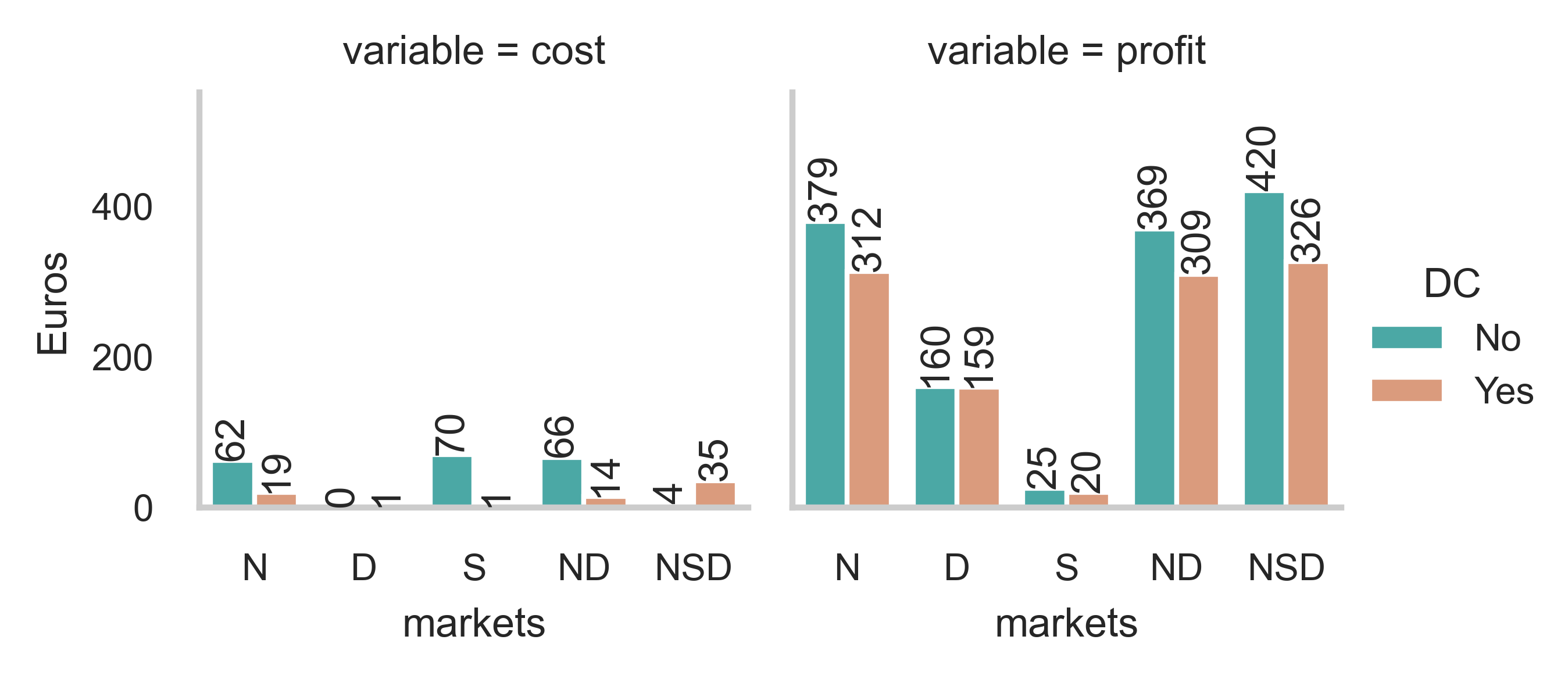}
    \caption{Costs and earnings for single-market participation \label{fig:opt_r14}}
\end{figure}

\subsection{BESS Services in Multiple Markets}

For $m = ND$, in both DC and no-DC case, bids are submitted and accepted 92\% of the time. However, the markets in which the bids are submitted differ. In no-DC case, the BESS prefers FCR-N market over FCR-D market. Thus, for 92\% of the time $m = N$, 8\% of the time $m = 0$, and 0\% of the time $m = D$. Contrarily, in DC case, BESS prefers FCR-D market over FCR-N market. Thus, for 25\% of the time $m = N$, 8\% of the time $m = 0$, and 67\% of the time $m = D$. This is because, as observed in the single-market participation, FCR-D has lower availability payments as compared to FCR-N and becomes attractive for a BESS only when its low energy demands resulting in lower degradation costs are accounted for. Overall profits of $m = N$ and $m = ND$ are similar.

Finally, when $m = NDS$, bids are respectively submitted and accepted 100\% and 79\% of the time for no-DC, and 92\% and 62\% for DC case. Moreover, in both no-DC and DC case, $m = \text{S-CH}$ is avoided. However, $m = \text{S-DCH}$ for 8\% of the time in no-DC case. This indicates that at certain hours, payments from $\text{S-DCH}$ are higher than $m = N$.

\begin{figure}[ht]
    \centering
    % Placeholder for costs and earnings figure
    \includegraphics[width=\linewidth]{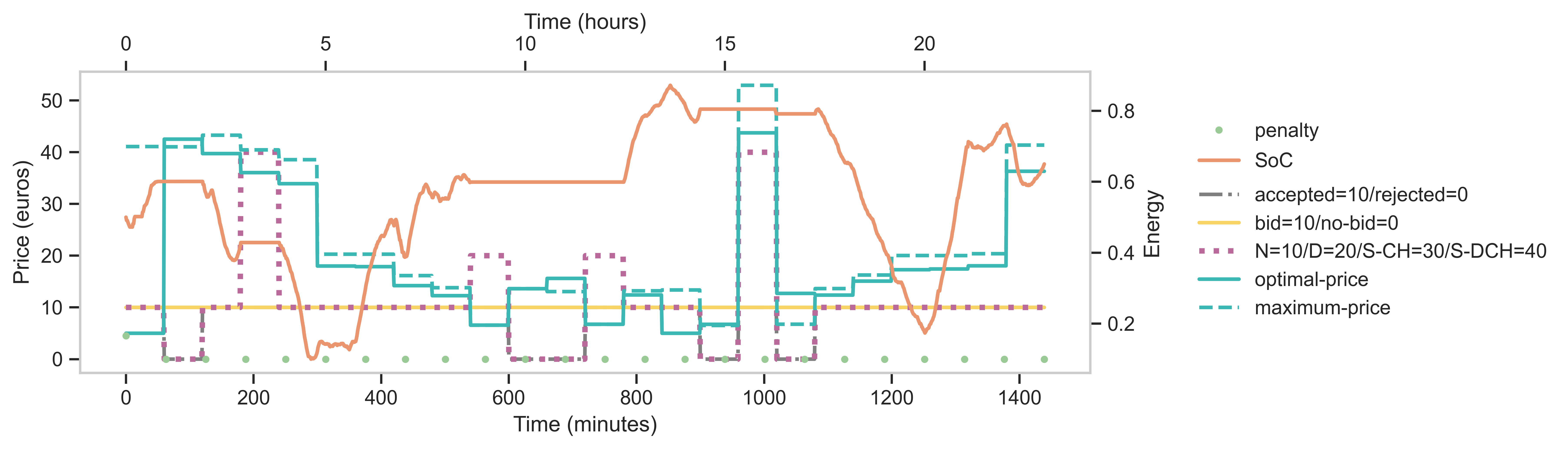}
    \caption{Charge/Discharge profile of BESS in multiple markets no-DC \label{fig:opt_r9}}
\end{figure}

\begin{figure}[ht]
    \centering
    % Placeholder for costs and earnings figure
    \includegraphics[width=\linewidth]{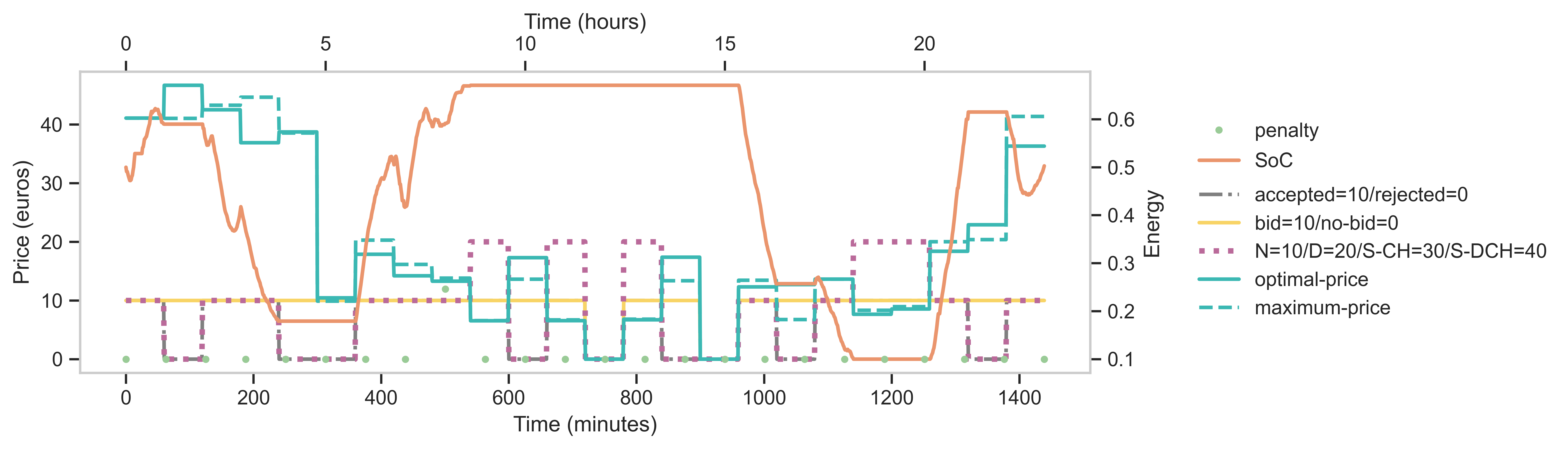}
    \caption{Charge/Discharge profile of BESS in multiple markets DC \label{fig:opt_r10}}
\end{figure}

\begin{figure}[ht]
    \centering
    % Placeholder for multi-market participation figure
    \includegraphics[width=1\linewidth]{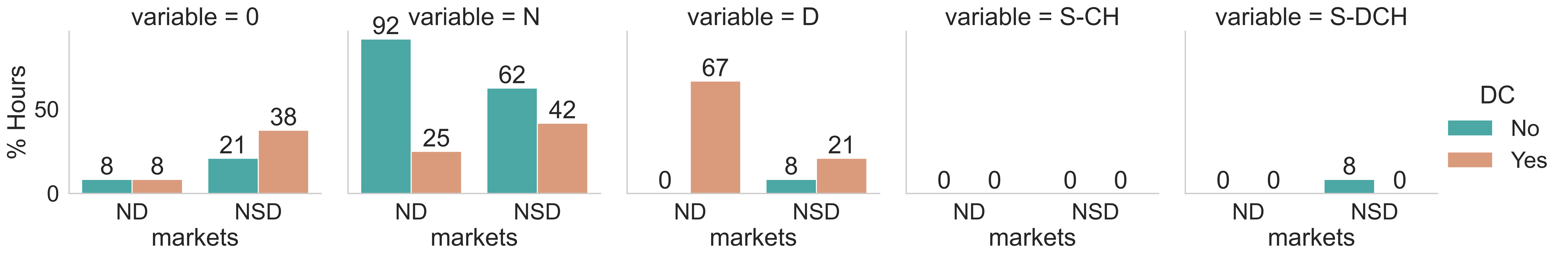}
    \caption{Percentage of hours in each market for stacked BESS services \label{fig:opt_r13}}
\end{figure}

\section{Conclusions}

In this paper we developed a stochastic optimization model to quantify the business value of BESS participation in different electricity markets. Our results showed that including the degradation costs in the optimization problem not only reduced the overall profits, but also the percentage of bid submission hours, both in single and multi-market BESS applications. However, reduction in profits was not linearly affected by the reduction in bid submission and acceptance hours. This implied that by identifying the suitable service switching hours and adequately limiting the bid submission hours can help to lower the impact of degradation costs on profits. However, this finding was mainly applicable to higher energy demand markets (HEDM) such as FCR-N and the spot market. For lower energy demand markets (LEDM), such as FCR-D, limiting bid submission hours was not desirable. The results thus imply that while bidding in frequency markets with 90\% power capacity, and spot markets with 40\% energy capacity, the frequency markets are a more attractive business solution.


\begin{thebibliography}{00}

\bibitem{b1} IRENA, \emph{Global energy transformation: A roadmap to 2050 (2019 edition)}, 2019.

\bibitem{b2} M. Wråke \emph{et al.}, “NORDIC CLEAN ENERGY SCENARIOS: Solutions for Carbon Neutrality,” 2021, doi:10.6027/ner2021-01.

\bibitem{b3} W. Yang \emph{et al.}, “Burden on hydropower units for short-term balancing of renewable power systems,” \emph{Nat. Commun.}, vol.~9, no.~1, pp.~1–12, 2018, doi:10.1038/s41467-018-05060-4.

\bibitem{b4} J. Figgener \emph{et al.}, “The development of stationary battery storage systems in Germany – A market review,” \emph{J. Energy Storage}, vol.~29, p.~101153, 2020, doi:10.1016/j.est.2019.101153.

\bibitem{b5} J. Badeda, J. Meyer, and D. U. Sauer, “Modeling the influence of installed battery energy storage systems on the German frequency containment reserve market,” in \emph{NEIS 2017 – Conf. Sustain. Energy Supply Energy Storage Syst.}, 2020, pp.~315–321.

\bibitem{b6} Z. Hameed, S. Hashemi, and C. Træholt, “Applications of AI-Based Forecasts in Renewable Based Electricity Balancing Markets,” in \emph{Proc. 2021 IEEE ICIT}, pp.~579–584, 2021, doi:10.1109/icit46573.2021.9453469.

\bibitem{b7} M. G. Pollitt and K. L. Anaya, “Competition in Markets for Ancillary Services? The Implications of Rising Distributed Generation,” \emph{Energy J.}, vol.~41, no.~01, 2020, doi:10.5547/01956574.42.si1.mpol.


\bibitem{b8} Asian Development Bank, \emph{Handbook on Battery Energy Storage System}, 2018.

\bibitem{b9} E. Namor and F. Sossan, “Control of Battery Storage Systems for the Simultaneous Provision of Multiple Services,” \emph{IEEE Trans. Smart Grid}, vol.~10, no.~3, pp.~2799–2808, 2020, doi:10.1109/TSG.2018.2810781.

\bibitem{b10} Y. Shi, B. Xu, D. Wang, and B. Zhang, “Using Battery Storage for Peak Shaving and Frequency Regulation: Joint Optimization for Superlinear Gains,” in \emph{2018 IEEE PESGM}, 2018, doi:10.1109/PESGM.2018.8586227.

\bibitem{b11} F. Gerini \emph{et al.}, “Optimal grid-forming control of battery energy storage systems providing multiple services: Modeling and experimental validation,” \emph{Electr. Power Syst. Res.}, vol.~212, p.~108567, 2022, doi:10.1016/j.epsr.2022.108567.

\bibitem{b12} H. Zhao, M. Hong, W. Lin, and K. A. Loparo, “Voltage and frequency regulation of microgrid with battery energy storage systems,” \emph{IEEE Trans. Smart Grid}, vol.~10, no.~1, pp.~414–424, 2019, doi:10.1109/TSG.2017.2741668.

\bibitem{b13} Z. Hameed, C. Træholt, and S. Hashemi, “Investigating the participation of battery energy storage systems in the Nordic ancillary services markets from a business perspective,” \emph{J. Energy Storage}, vol.~58, p.~106464, 2023, doi:10.1016/j.est.2022.106464.

\bibitem{b14} Z. Hameed, S. Hashemi, H. H. Ipsen, and C. Træholt, “A business-oriented approach for battery energy storage placement in power systems,” \emph{Appl. Energy}, vol.~298, p.~117186, 2021, doi:10.1016/j.apenergy.2021.117186.

\bibitem{b15} Z. Hameed, S. Hashemi, and C. Træholt, “Site selection criteria for battery energy storage in power systems,” in \textit{Proc. 2020 IEEE Canadian Conf. on Electrical and Computer Engineering (CCECE)}, 2020, pp. 1–5.

\bibitem{b16} Z. Hameed, S. Hashemi, H. H. Ipsen, and C. Træholt, “Placement of Battery Energy Storage for Provision of Grid Services–A Bornholm Case Study,” in \textit{Proc. 2021 IEEE 9th Int. Conf. on Smart Energy Grid Engineering (SEGE)}, 2021, pp. 1–6.

\bibitem{b17} Z. Hameed, M. Pollitt, P. Kattuman, and C. Træholt, “Frequency markets and the problem of pre-dictability,” Tech. Rep., Faculty of Economics, University of Cambridge, 2023.

\bibitem{b18} Y. Zhang, Y. Xu, H. Yang, Z. Y. Dong, and R. Zhang, “Optimal Whole-Life-Cycle Planning of Battery Energy Storage for Multi-Functional Services in Power Systems,” \emph{IEEE Trans. Sustain. Energy}, vol.~11, no.~4, pp.~2077–2086, 2020, doi:10.1109/TSTE.2019.2942066.

\bibitem{b19} T. Sayfutdinov and P. Vorobev, “Optimal utilization strategy of the LiFePO4 battery storage,” \emph{Appl. Energy}, vol.~316, 2022, doi:10.1016/j.apenergy.2022.119080.


\end{thebibliography}
\end{document}